\documentclass[10pt]{article}
\usepackage{cite}
\usepackage{mathrsfs}
\usepackage{amsfonts}
\usepackage{amsmath}
\usepackage{amsfonts,amssymb,color}
\usepackage{dsfont}
\usepackage{curves}
\usepackage{mathrsfs}
\usepackage{pifont}
\usepackage{amssymb}
\allowdisplaybreaks

\numberwithin{equation}{section}

\date{}

\def\BigRoman{\uppercase\expandafter{\romannumeral\number\count 255 }}
\def\Romannumeral{\afterassignment\BigRoman\count255=}

\begin{document}
\title{Sufficient conditions for spanning trees with bounded total excess in tough graphs
}
\author{\small Sizhong Zhou$^{1}$\footnote{Corresponding author. E-mail address: zsz\_cumt@163.com}, Chunjie Sun$^{2}$, Zhiren Sun$^{2,3}$\\
\small $1$. School of Science, Jiangsu University of Science and Technology,\\
\small Zhenjiang, Jiangsu 212100, China\\
\small $2$. Basic Teaching Department, Xinjiang University of Political Science and Law,\\
\small Tumushuke, Xinjiang 843900, China\\
\small $3$. School of Mathematical Sciences, Nanjing Normal University,\\
\small Nanjing, Jiangsu 210023, China\\
}

\maketitle
\begin{abstract}
\noindent Let $G$ be a connected graph and $T$ a spanning tree of $G$. The toughness of $G$ is defined by $t(G)=\min\Big\{\frac{|S|}{c(G-S)}:S\subseteq V(G), c(G-S)\geq2\Big\}$ for $G\neq K_n$. A graph
$G$ is said to be $t$-tough if $t(G)\geq t$. The size, the adjacency spectral radius and the signless Laplacian spectral radius of $G$ are denoted by $e(G)$, $\rho(G)$ and $q(G)$, respectively. The
$k$-excess of a vertex $v$ in $T$ is defined as $\max\{0,d_T(v)-k\}$. The total $k$-excess $\mbox{te}(T,k)$ is defined by $\mbox{te}(T,k)=\sum\limits_{v\in V(T)}{\max\{0,d_T(v)-k\}}$. A tree $T$ is
called a $k$-tree if $d_T(v)\leq k$ for every $v\in V(T)$, namely, the maximum degree of a $k$-tree is at most $k$. It is obvious that $T$ is a spanning $k$-tree if and only if $\mbox{te}(T,k)=0$. This
paper investigates a generalization of spanning $k$-trees using a concept called total $k$-excess and puts forward three lower bounds based on $e(G)$, $\rho(G)$ and $q(G)$ in a connected
$\frac{k}{k(k-2)+b+3}$-tough graph $G$ to guarantee that $G$ has a spanning tree $T$ with $\mbox{te}(T,k)\leq b$, where $k\geq3$ and $b\geq0$ are two integers. As corollaries, we obtain three sufficient
conditions for a connected $\frac{k}{k(k-2)+3}$-tough graph to has a spanning $k$-tree.
\\
\begin{flushleft}

(2020) Mathematics Subject Classification: 05C50, 05C05
\end{flushleft}
\end{abstract}

\section{Introduction}

Throughout this paper, we only deal with finite, undirected and simple graphs. Let $G$ be a graph. The vertex set and the edge set of $G$ are denoted by $V(G)$ and $E(G)$, respectively. We use $e(G)$
and $n$ to denote the size and the order of $G$, respectively. Then $|E(G)|=e(G)$ and $|V(G)|=n$. Let $d_G(v)$ and $c(G)$ denote the degree of a vertex $v$ in $G$ and the number of components in $G$,
respectively. Given a subset $S\subseteq V(G)$, we denote by $G[S]$ the subgraph of $G$ induced by $S$, and by $G-S$ the graph formed from $G$ by deleting the vertices in $S$ together with the edges
incident with vertices in $S$. Let $S$ and $T$ be two disjoint vertex subset of $G$. We denote by $e_G(S,T)$ the number of edges with one end in $S$ and the other end in $T$. The complete graph with
$n$ vertices is denoted by $K_n$. Given two graphs $G_1$ and $G_2$, the union of $G_1$ and $G_2$, denoted by $G_1\cup G_2$, is the graph with vertex set $V(G_1)\cup V(G_2)$ and edge set
$E(G_1)\cup E(G_2)$. The join of $G_1$ and $G_2$, denoted by $G_1\vee G_2$, is obtained from $G_1\cup G_2$ by adding all possible edges between $V(G_1)$ and $V(G_2)$. The union of $k$ copies of $G$ is
denoted by $kG$, where $k\geq1$ is an integer. The toughness of a graph $G$ is defined by
$$
t(G)=\min\Big\{\frac{|S|}{c(G-S)}:S\subseteq V(G), c(G-S)\geq2\Big\}
$$
for $G\neq K_n$. A graph $G$ is said to be $t$-tough if $t(G)\geq t$.

Let $G$ be a connected graph with $V(G)=\{v_1,v_2,\ldots,v_n\}$. The adjacency matrix of $G$, denoted by $A(G)$, is defined as $A(G)=(a_{ij})_{n\times n}$, where $a_{ij}=1$ if $v_i$ and $v_j$ are
adjacent in $G$, $a_{ij}=0$ otherwise. The largest eigenvalue of $A(G)$ is called the adjacency spectral radius of $G$, and written as $\rho(G)$. The signless Laplacian matrix of $G$ is defined as
$Q(G)=D(G)+A(G)$, where $D(G)=\mbox{diag}(d_G(v_1),d_G(v_2),\ldots,d_G(v_n))$. The signless Laplacian spectral radius of $G$, denoted by $q(G)$, is the largest eigenvalue of $Q(G)$. The spectral radii
of graphs have been extensively investigated, and some properties on the adjacency spectral radius \cite{B,Os,Ws2,Zs1,Zs2,Zs3} and the signless Laplacian spectral radius \cite{PK,Ws1,ZL,ZLB,ZZZL} of
$G$ were obtained.

Given a connected graph $G$ and an integer $k\geq2$, a spanning $k$-tree $T$ of $G$ is a spanning tree in which each vertex has degree at most $k$. Obviously, a Hamiltonian path is a spanning 2-tree.
And so, a spanning $k$-tree of $G$ is a natural generalization of a Hamiltonian path. The $k$-excess of a vertex $v$ in $T$ is defined by $\max\{0,d_T(v)-k\}$. We write $\mbox{te}(T,k)=\sum\limits_{v\in V(T)}\max\{0,d_T(v)-k\}$,
which is called the total $k$-excess of $T$. In fact, $T$ is a spanning $k$-tree if and only if $\mbox{te}(T,k)=0$. The study on the existence of a spanning tree with bounded total excess in a graph has
an important theoretical significance in graph theory, as it naturally generalizes the classical Hamiltonian path problem and bridges the gap between global degree distributions and structural constraints.

Kaneko, Kano and Suzuki \cite{KKS} investigated the existence of a spanning tree with leaf distance at least 4 in a graph. Wu \cite{Wc} provided two sufficient conditions to guarantee that a graph admits
a spanning tree with leaf degree at most $k$ with respect to its size and adjacency spectral radius. Wu, Zhou and Liu \cite{WZL} created a lower bound on the signless Laplacian spectral radius to ensure
that a bipartite graph $G$ with bipartition $(A,B)$ contains a spanning tree $T$ with $d_T(v)\geq k$ for each $v\in A$. Ota and Sugiyama \cite{OS} obtained a sufficient condition for a graph to possess
a spanning $k$-tree based on forbidden subgraphs. Furuya et al. \cite{FMMMTY} presented some degree sum conditions for a star-free graph to contain a spanning $k$-tree. Rose \cite{R} characterised a
graph with a spanning $k$-tree. Gu and Liu \cite{GL} gave a Laplacian eigenvalue condition for the existence of a spanning $k$-tree in a graph. Fan et al. \cite{FGHL} provided a characterization of a
graph with a spanning $k$-tree in terms of its adjacency spectral radius and signless Laplacian spectral radius. Zhou and Wu \cite{ZW} established an upper bound for the distance spectral radius to guarantee
that a graph has a spanning $k$-tree. Wu and Zhou \cite{WZ} created a connection between the adjacency spectral radius and spanning $k$-trees in binding graphs. Liu, Fan and Shu \cite{LFS} showed two
sufficient conditions based on the adjacency spectral radius and the signless Laplacian spectral radius for a tough graph to have a spanning $k$-tree. Ohnishi, Ota and Ozeki \cite{OOO} obtained some
sufficient conditions for the existence of a spanning tree with bounded total $k$-excess in a graph. Maezawa, Tsugaki and Yashima \cite{MTY} posed a new closure concept and applied it to investigate
spanning trees with bounded total $k$-excess of graphs. Ozeki \cite{Oa} proposed a toughness condition for the existence of a spanning tree with bounded total $k$-excess in a graph.

Win \cite{Win} obtained the following result, which proposes a sufficient condition for the existence of a spanning $k$-tree in a graph. A short proof was provided by Ellingham and Zha \cite{EZ}.

\medskip

\noindent{\textbf{Theorem 1.1}} (Win \cite{Win}, Ellingham and Zha \cite{EZ}). Let $k$ be an integer with $k\geq2$, and let $G$ be a connected graph. If
$$
c(G-S)\leq(k-2)|S|+2
$$
for every subset $S\subseteq V(G)$, then $G$ contains a spanning $k$-tree.

\medskip

Enomoto, Ohnishi and Ota \cite{EOO} presented a generalization of Theorem 1.1 with respect to bounded total $k$-excess $\mbox{te}(T,k)$.

\medskip

\noindent{\textbf{Theorem 1.2}} (Enomoto, Ohnishi and Ota \cite{EOO}). Let $b$ and $k$ be two nonnegative integers with $k\geq2$, and let $G$ be a connected graph. If
$$
c(G-S)\leq(k-2)|S|+b+2
$$
for each subset $S\subseteq V(G)$, then $G$ contains a spanning tree $T$ with $\mbox{te}(T,k)\leq b$.

\medskip

Zhou \cite{Za} provided an adjacency spectral condition for a graph to possess a spanning tree with bounded total $k$-excess.

\medskip

\noindent{\textbf{Theorem 1.3}} (Zhou \cite{Za}). Let $k$ and $b$ be two nonnegative integers with $k\geq\max\{5,b+3\}$ and $(b,k)\neq(2,5)$, and let $G$ be a connected graph of order $n\geq k+b+2$. If
$$
\rho(G)\geq\rho(K_1\vee(K_{n-k-b-1}\cup(k+b)K_1)),
$$
then $G$ has a spanning tree $T$ with $\mbox{te}(T,k)\leq b$, unless $G=K_1\vee(K_{n-k-b-1}\cup(k+b)K_1)$.

\medskip

Gong and Zhang \cite{GZ} arose two sufficient conditions, with respect to the size and the signless Laplacian spectral radius, for a graph to have a spanning tree with bounded total $k$-excess.

\medskip

\noindent{\textbf{Theorem 1.4}} (Gong and Zhang \cite{GZ}). Let $k$ and $b$ be two nonnegative integers with $k>\max\{4,\frac{5+\sqrt{8b+17}}{2}\}$. If $G$ is a connected graph of order $n\geq k+b+2$
such that
$$
e(G)\geq\binom{n-k-b}{2}+k+b,
$$
then $G$ has a spanning tree $T$ with $\mbox{te}(T,k)\leq b$, unless $G=K_1\vee(K_{n-k-b-1}\cup(k+b)K_1)$.

\medskip

\noindent{\textbf{Theorem 1.5}} (Gong and Zhang \cite{GZ}). Let $k$ and $b$ be two nonnegative integers with $k\geq\max\{4,b+3\}$. If $G$ is a connected graph of order $n\geq2k+2b+10$ such that
$$
q(G)\geq q(K_1\vee(K_{n-k-b-1}\cup(k+b)K_1)),
$$
then $G$ has a spanning tree $T$ with $\mbox{te}(T,k)\leq b$, unless $G=K_1\vee(K_{n-k-b-1}\cup(k+b)K_1)$.

\medskip

Motivated by Theorems 1.2--1.5, we create some sufficient conditions based on the size, the adjacency spectral radius and the signless Laplacian spectral radius to ensure the existence of a
spanning tree $T$ with bounded total $k$-excess in a connected $\frac{k}{k(k-2)+b+3}$-tough graph. We are to show our main results as follows.

\medskip

\noindent{\textbf{Theorem 1.6.}} Let $k$ and $b$ be two nonnegative integers with $k\geq3$, and let $G$ be a connected $\frac{k}{k(k-2)+b+3}$-tough graph of order $n\geq k^{2}+2k+b+4+\frac{4b+9}{k-2}$.
If
$$
e(G)\geq e(K_k\vee(K_{n-k(k-1)-b-2}\cup(k(k-2)+b+2)K_1)),
$$
then $G$ has a spanning tree $T$ with $\mbox{te}(T,k)\leq b$, unless $G=K_k\vee(K_{n-k(k-1)-b-2}\cup(k(k-2)+b+2)K_1)$.

\medskip

\noindent{\textbf{Theorem 1.7.}} Let $k\geq3$ and $b\geq0$ be two integers, and let $G$ be a connected $\frac{k}{k(k-2)+b+3}$-tough graph of order $n\geq2k^{2}+k+3b+3$. If
$$
\rho(G)\geq\rho(K_k\vee(K_{n-k(k-1)-b-2}\cup(k(k-2)+b+2)K_1)),
$$
then $G$ has a spanning tree $T$ with $\mbox{te}(T,k)\leq b$, unless $G=K_k\vee(K_{n-k(k-1)-b-2}\cup(k(k-2)+b+2)K_1)$.

\medskip

\noindent{\textbf{Theorem 1.8.}} Let $k\geq3$ and $b\geq0$ be two integers, and let $G$ be a connected $\frac{k}{k(k-2)+b+3}$-tough graph of order $n\geq\frac{k^{4}+2bk^{2}+k+b^{2}+b-4}{2(k-2)}$. If
$$
q(G)\geq q(K_k\vee(K_{n-k(k-1)-b-2}\cup(k(k-2)+b+2)K_1)),
$$
then $G$ has a spanning tree $T$ with $\mbox{te}(T,k)\leq b$, unless $G=K_k\vee(K_{n-k(k-1)-b-2}\cup(k(k-2)+b+2)K_1)$.

\medskip

We immediately obtain the following corollaries by setting $b=0$ in Theorems 1.6--1.8.

\medskip

\noindent{\textbf{Corollary 1.9.}} Let $k$ be an integer with $k\geq3$, and let $G$ be a connected $\frac{k}{k(k-2)+3}$-tough graph of order $n\geq k^{2}+2k+4+\frac{9}{k-2}$. If
$$
e(G)\geq e(K_k\vee(K_{n-k(k-1)-2}\cup(k(k-2)+2)K_1)),
$$
then $G$ has a spanning $k$-tree, unless $G=K_k\vee(K_{n-k(k-1)-2}\cup(k(k-2)+2)K_1)$.

\medskip

\noindent{\textbf{Corollary 1.10.}} Let $k\geq3$ be an integer, and let $G$ be a connected $\frac{k}{k(k-2)+3}$-tough graph of order $n\geq2k^{2}+k+3$. If
$$
\rho(G)\geq\rho(K_k\vee(K_{n-k(k-1)-2}\cup(k(k-2)+2)K_1)),
$$
then $G$ has a spanning $k$-tree, unless $G=K_k\vee(K_{n-k(k-1)-2}\cup(k(k-2)+2)K_1)$.

\medskip

\noindent{\textbf{Corollary 1.11.}} Let $k\geq3$ be an integer, and let $G$ be a connected $\frac{k}{k(k-2)+3}$-tough graph of order $n\geq\frac{k^{4}+k-4}{2(k-2)}$. If
$$
q(G)\geq q(K_k\vee(K_{n-k(k-1)-2}\cup(k(k-2)+2)K_1)),
$$
then $G$ has a spanning $k$-tree, unless $G=K_k\vee(K_{n-k(k-1)-2}\cup(k(k-2)+2)K_1)$.

\medskip

\section{Preliminary lemmas}

In this section, we list some lemmas which will be useful in the proofs of our main results.

\medskip

\noindent{\textbf{Lemma 2.1}} (Zheng et al. \cite{ZLLW}). Let $\sum\limits_{i=1}^{t}n_i+s=n$. If $n_1\geq n_2\geq\cdots\geq n_t\geq p\geq1$, then
$$
e(K_s\vee(K_{n_1}\cup K_{n_2}\cup\cdots\cup K_{n_t}))\leq e(K_s\vee(K_{n-s-p(t-1)}\cup(t-1)K_p)),
$$
with equality occurring if and only if $(n_1,n_2,\ldots,n_t)=(n-s-p(t-1),p,\ldots,p)$.

\medskip

\noindent{\textbf{Lemma 2.2}} (Li and Feng \cite{LF}). Let $H$ be a subgraph of a connected graph $G$. Then
$$
\rho(G)\geq\rho(H),
$$
where the equality follows if and only if $G=H$.

\medskip

\noindent{\textbf{Lemma 2.3}} (Fan et al.\cite{FGHL}). Let $\sum\limits_{i=1}^{t}n_i+s=n$. If $n_1\geq n_2\geq\cdots\geq n_t\geq1$, then
$$
\rho(K_s\vee(K_{n_1}\cup K_{n_2}\cup\cdots\cup K_{n_t}))\leq\rho(K_s\vee(K_{n-s-t+1}\cup(t-1)K_1)),
$$
with equality holding if and only if $(n_1,n_2,\ldots,n_t)=(n-s-t+1,1,\ldots,1)$.

\medskip

\medskip

\noindent{\textbf{Lemma 2.4}} (Shen et al. \cite{SYZL}). If $H$ is a subgraph of a connected graph $G$, then
$$
q(G)\geq q(H),
$$
with equality if and only if $G=H$.

\medskip

\noindent{\textbf{Lemma 2.5}} (Zheng et al. \cite{ZLLW}). Let $\sum\limits_{i=1}^{t}n_i+s=n$. If $n_1\geq n_2\geq\cdots\geq n_t\geq p\geq1$, then
$$
q(K_s\vee(K_{n_1}\cup K_{n_2}\cup\cdots\cup K_{n_t}))\leq q(K_s\vee(K_{n-s-p(t-1)}\cup(t-1)K_p)),
$$
with equality occurring if and only if $(n_1,n_2,\ldots,n_t)=(n-s-p(t-1),p,\ldots,p)$.

\medskip

\noindent{\textbf{Lemma 2.6}} (Das \cite{Das}). Let $G$ be a graph of order $n$. Then
$$
q(G)\leq\frac{2e(G)}{n-1}+n-2.
$$

\medskip

Let $M$ be a real matrix whose rows and columns are indexed by $\mathcal{N}=\{1,2,\ldots,n\}$. For a given partition $\pi:\mathcal{N}=\mathcal{N}_1\cup\mathcal{N}_2\cup\cdots\cup\mathcal{N}_r$, the matrix $M$
can be correspondingly partitioned into
\begin{align*}
M=\left(
  \begin{array}{cccc}
    M_{11} & M_{12} & \cdots & M_{1r}\\
    M_{21} & M_{22} & \cdots & M_{2r}\\
    \vdots & \vdots & \ddots & \vdots\\
    M_{r1} & M_{r2} & \cdots & M_{rr}\\
  \end{array}
\right).
\end{align*}
The quotient matrix of $M$ in terms of $\pi$ is defined by $M_{\pi}=(m_{ij})_{r\times r}$, where $m_{ij}$ denotes the average row sum of the block $M_{ij}$ for each $1\leq i,j\leq r$. The partition $\pi$ is
called equitable if the row sum of every block $M_{ij}$ of $M$ is constant for $1\leq i,j\leq r$.

\medskip

\noindent{\textbf{Lemma 2.7}} (You et al. \cite{YYSX}). Let $M$ be a real $n\times n$ matrix with an equitable partition $\pi$, and let $M_{\pi}$ be the corresponding quotient matrix. Then each
eigenvalue of $M_{\pi}$ is an eigenvalue of $M$. Furthermore, if $M$ is nonnegative and irreducible, then the largest eigenvalues of $M$ and $M_{\pi}$ are equal.

\medskip

\section{The proof of Theorem 1.6}

\noindent{\it Proof of Theorem 1.6.} Suppose to the contrary that a connected $\frac{k}{k(k-2)+b+3}$-tough graph $G$ has no spanning tree $T$ with $\mbox{te}(T,k)\leq b$. By Theorem 1.2, there exists some
vertex subset $S$ of $G$ satisfying $c(G-S)\geq(k-2)|S|+b+3$. Thus, we deduce $n\geq c(G-S)+|S|\geq(k-2)|S|+b+3+|S|=(k-1)|S|+b+3$. Notice that $G$ is connected. Then we obtain $S\neq\emptyset$. For convenience,
we write $|S|=s$ and $c(G-S)=c$. Then $c\geq(k-2)s+b+3$, $n\geq(k-1)s+b+3$ and $G$ is a spanning subgraph of $G_1=K_s\vee(K_{n_1}\cup K_{n_2}\cup\cdots\cup K_{n_{(k-2)s+b+3}})$ for some positive integers
$n_1,n_2,\ldots,n_{(k-2)s+b+3}$ with $\sum\limits_{i=1}^{(k-2)s+b+3}n_i+s=n$. Without loss of generality, we assume that $n_1\geq n_2\geq\cdots\geq n_{(k-2)s+b+3}\geq1$. Hence, we conclude
\begin{align}\label{eq:3.1}
e(G)\leq e(G_1),
\end{align}
where the equality follows if and only if $G=G_1$.

Let $G_2=K_s\vee(K_{n-(k-1)s-b-2}\cup((k-2)s+b+2)K_1)$. According to Lemma 2.1, we possess
\begin{align}\label{eq:3.2}
e(G_1)\leq e(G_2),
\end{align}
with equality if and only if $(n_1,n_2\ldots,n_{(k-2)s+b+3})=(n-(k-1)s-b-2,1\ldots,1)$.

Recall that $|S|=s$ and $c(G-S)\geq(k-2)|S|+b+3$. In terms of $k\geq3$, $b\geq0$ and the definition of $\frac{k}{k(k-2)+b+3}$-tough graphs, we have
$$
\frac{k}{k(k-2)+b+3}\leq t(G)\leq\frac{|S|}{c(G-S)}\leq\frac{|S|}{(k-2)|S|+b+3},
$$
which yields that $s=|S|\geq k$.

If $s=k$, then we have $G_2=K_k\vee(K_{n-k(k-1)-b-2}\cup(k(k-2)+b+2)K_1)$, and so $e(G_2)=e(K_k\vee(K_{n-k(k-1)-b-2}\cup(k(k-2)+b+2)K_1))$. Together with \eqref{eq:3.1} and \eqref{eq:3.2}, we obtain
$$
e(G)\leq e(K_k\vee(K_{n-k(k-1)-b-2}\cup(k(k-2)+b+2)K_1)),
$$
with equality following if and only if $G=K_k\vee(K_{n-k(k-1)-b-2}\cup(k(k-2)+b+2)K_1)$. Combining this with the condition $e(G)\geq e(K_k\vee(K_{n-k(k-1)-b-2}\cup(k(k-2)+b+2)K_1))$ of the theorem,
we conclude $G=K_k\vee(K_{n-k(k-1)-b-2}\cup(k(k-2)+b+2)K_1)$, which contradicts that $G\neq K_k\vee(K_{n-k(k-1)-b-2}\cup(k(k-2)+b+2)K_1)$. Next, we are to consider $s\geq k+1$.

Let $G_*=K_k\vee(K_{n-k(k-1)-b-2}\cup(k(k-2)+b+2)K_1)$. Recall that $G_2=K_s\vee(K_{n-(k-1)s-b-2}\cup((k-2)s+b+2)K_1)$ and $n\geq(k-1)s+b+3$. By a simple computation, we obtain
\begin{align*}
e(G_2)-e(G_*)=&\binom{n-(k-2)s-b-2}{2}+s\Big((k-2)s+b+2\Big)\\
&-\binom{n-k(k-2)-b-2}{2}-k\Big(k(k-2)+b+2\Big)\\
=&\frac{(n-(k-2)s-b-2)(n-(k-2)s-b-3)}{2}+s\Big((k-2)s+b+2\Big)\\
&-\frac{(n-k(k-2)-b-2)(n-k(k-2)-b-3)}{2}-k\Big(k(k-2)+b+2\Big)\\
=&\frac{1}{4}(s-k)\Big(-4(k-2)n+2k(k-2)s+2k^{3}-4k^{2}+4bk+10k-4b-12\Big)\\
\leq&\frac{1}{4}(s-k)\Big(-(k-2)n-3(k-2)((k-1)s+b+3)+2k(k-2)s\\
&+2k^{3}-4k^{2}+4bk+10k-4b-12\Big)\\
& \ \ \ \ \ \ \ \ \ \ \ \ (\mbox{since} \ k\geq3, \ s\geq k+1 \ \mbox{and} \ n\geq(k-1)s+b+3)\\
=&\frac{1}{4}(s-k)\Big(-(k-2)n-(k-2)(k-3)s+2k^{3}-4k^{2}+bk+k+2b+6\Big)\\
\leq&\frac{1}{4}(s-k)\Big(-(k-2)n-(k-2)(k-3)(k+1)+2k^{3}-4k^{2}+bk+k+2b+6\Big)\\
& \ \ \ \ \ \ \ \ \ \ \ \ (\mbox{since} \ k\geq3 \ \mbox{and} \ s\geq k+1)\\
=&\frac{1}{4}(s-k)\Big(-(k-2)n+k^{3}+bk+2b\Big)\\
<&0 \ \ \ \ \ \ \ \ \ \ \ \ \Big(\mbox{since} \ k\geq3, \ s\geq k+1 \ \mbox{and} \ n\geq k^{2}+2k+b+4+\frac{4b+9}{k-2}\Big),
\end{align*}
which implies
\begin{align}\label{eq:3.3}
e(G_2)<e(G_*)=e(K_k\vee(K_{n-k(k-1)-b-2}\cup(k(k-2)+b+2)K_1)).
\end{align}
Applying \eqref{eq:3.1}, \eqref{eq:3.2} and \eqref{eq:3.3}, we infer
$$
e(G)\leq e(G_1)\leq e(G_2)<e(K_k\vee(K_{n-k(k-1)-b-2}\cup(k(k-2)+b+2)K_1)),
$$
which is a contradiction to $e(G)\geq e(K_k\vee(K_{n-k(k-1)-b-2}\cup(k(k-2)+b+2)K_1))$. Theorem 1.6 is proved. \hfill $\Box$

\medskip

\section{The proof of Theorem 1.7}

\noindent{\it Proof of Theorem 1.7.} Suppose to the contrary that a connected $\frac{k}{k(k-2)+b+3}$-tough graph $G$ has no spanning tree $T$ with $\mbox{te}(T,k)\leq b$. By Theorem 1.2, there exists some
vertex subset $S$ of $G$ satisfying $c(G-S)\geq(k-2)|S|+b+3$. Thus, we obtain $n\geq c(G-S)+|S|\geq(k-2)|S|+b+3+|S|=(k-1)|S|+b+3$. Since $G$ is connected, we deduce $S\neq\emptyset$. For convenience,
we write $|S|=s$ and $c(G-S)=c$. Then $c\geq(k-2)s+b+3$, $n\geq(k-1)s+b+3$ and $G$ is a spanning subgraph of $G_1=K_s\vee(K_{n_1}\cup K_{n_2}\cup\cdots\cup K_{n_{(k-2)s+b+3}})$ for some positive integers
$n_1,n_2,\ldots,n_{(k-2)s+b+3}$ with $\sum\limits_{i=1}^{(k-2)s+b+3}n_i+s=n$. Without loss of generality, we assume that $n_1\geq n_2\geq\cdots\geq n_{(k-2)s+b+3}\geq1$. In terms of Lemma 2.2, we get
\begin{align}\label{eq:4.1}
\rho(G)\leq\rho(G_1),
\end{align}
with equality if and only if $G=G_1$.

Let $G_2=K_s\vee(K_{n-(k-1)s-b-2}\cup((k-2)s+b+2)K_1)$. According to Lemma 2.3, we conclude
\begin{align}\label{eq:4.2}
\rho(G_1)\leq\rho(G_2),
\end{align}
where the equality occurs if and only if $(n_1,n_2\ldots,n_{(k-2)s+b+3})=(n-(k-1)s-b-2,1\ldots,1)$.

Recall that $|S|=s$ and $c(G-S)\geq(k-2)|S|+b+3$. According to $k\geq3$, $b\geq0$ and the definition of $\frac{k}{k(k-2)+b+3}$-tough graphs, we obtain
$$
\frac{k}{k(k-2)+b+3}\leq t(G)\leq\frac{|S|}{c(G-S)}\leq\frac{|S|}{(k-2)|S|+b+3},
$$
which implies $s=|S|\geq k$.

If $s=k$, then we obtain $G_2=K_k\vee(K_{n-k(k-1)-b-2}\cup(k(k-2)+b+2)K_1)$, and so $\rho(G_2)=\rho(K_k\vee(K_{n-k(k-1)-b-2}\cup(k(k-2)+b+2)K_1))$. Then it follows from \eqref{eq:4.1} and \eqref{eq:4.2}
that
$$
\rho(G)\leq\rho(K_k\vee(K_{n-k(k-1)-b-2}\cup(k(k-2)+b+2)K_1)),
$$
where the equality occurs if and only if $G=K_k\vee(K_{n-k(k-1)-b-2}\cup(k(k-2)+b+2)K_1)$. Together with the condition $\rho(G)\geq\rho(K_k\vee(K_{n-k(k-1)-b-2}\cup(k(k-2)+b+2)K_1))$ of the theorem, we
possess $G=K_k\vee(K_{n-k(k-1)-b-2}\cup(k(k-2)+b+2)K_1)$, which is a contradiction to $G\neq K_k\vee(K_{n-k(k-1)-b-2}\cup(k(k-2)+b+2)K_1)$. In what follows, we are to consider $s\geq k+1$.

Recall that $G_2=K_s\vee(K_{n-(k-1)s-b-2}\cup((k-2)s+b+2)K_1)$. The quotient matrix of $A(G_2)$ based on the partition $V(G_2)=V(K_{n-(k-1)s-b-2})\cup V(((k-2)s+b+2)K_1)\cup V(K_s)$ is given by
\begin{align*}
B_2=\left(
  \begin{array}{ccc}
    n-(k-1)s-b-3 & 0 & s\\
    0 & 0 & s\\
    n-(k-1)s-b-2 & (k-2)s+b+2 & s-1\\
  \end{array}
\right),
\end{align*}
and so its characteristic polynomial is written as
\begin{align*}
\varphi_{B_2}(x)=&x^{3}-(n-(k-2)s-b-4)x^{2}-(n+(k-2)s^{2}+(b-k+4)s-b-3)x\\
&+(k-2)s^{2}n+(b+2)sn-(k^{2}-3k+2)s^{3}-(2bk+5k-3b-8)s^{2}\\
&-(b^{2}+5b+6)s.
\end{align*}

Let $G_*=K_k\vee(K_{n-k(k-1)-b-2}\cup(k(k-2)+b+2)K_1)$. Then the quotient matrix of $A(G_*)$ according to the partition $V(G_*)=V(K_{n-k(k-1)-b-2})\cup V((k(k-2)+b+2)K_1)\cup V(K_k)$ can be written as
\begin{align*}
B_*=\left(
  \begin{array}{ccc}
    n-k(k-1)-b-3 & 0 & k\\
    0 & 0 & k\\
    n-k(k-1)-b-2 & k(k-2)+b+2 & k-1\\
  \end{array}
\right).
\end{align*}
The characteristic polynomial of $B_*$ equals
\begin{align*}
\varphi_{B_*}(x)=&x^{3}-(n-k^{2}+2k-b-4)x^{2}-(n+k^{3}-3k^{2}+(b+4)k-b-3)x\\
&+(k^{3}-2k^{2}+(b+2)k)n-k^{5}+3k^{4}-(2b+7)k^{3}+(3b+8)k^{2}\\
&-(b^{2}+5b+6)k.
\end{align*}

Notice that the two partitions $V(G_2)=V(K_{n-(k-1)s-b-2})\cup V(((k-2)s+b+2)K_1)\cup V(K_s)$ and $V(G_*)=V(K_{n-k(k-1)-b-2})\cup V((k(k-2)+b+2)K_1)\cup V(K_k)$ are both equitable. In terms of Lemma 2.7,
the largest root of $\varphi_{B_2}(x)=0$ is equal to $\rho(G_2)$ and the largest root of $\varphi_{B_*}(x)=0$ equals $\rho(G_*)$.

Since $G_*=K_k\vee(K_{n-k(k-1)-b-2}\cup(k(k-2)+b+2)K_1)$ contains $K_{n-k(k-2)-b-2}$ as its proper subgraph, it follows from Lemma 2.2 that
$$
\rho(G_*)>\rho(K_{n-k(k-2)-b-2})=n-k(k-2)-b-3.
$$
In order to verify $\rho(G_2)<\rho(G_*)$, it suffices to show $\varphi_{B_2}(x)-\varphi_{B_*}(x)>0$ for $x\geq n-k(k-2)-b-3$.

By a direct computation, we have
\begin{align}\label{eq:4.3}
\varphi_{B_2}(x)-\varphi_{B_*}(x)=(s-k)f_1(x),
\end{align}
where $f_1(x)=(k-2)x^{2}-((k-2)s+k^{2}-3k+b+4)x+(k-2)(s+k)n+(b+2)n-(k^{2}-3k+2)s^{2}-(k^{3}-3k^{2}+(2b+7)k-3b-8)s-k^{4}+3k^{3}-(2b+7)k^{2}+(3b+8)k-b^{2}-5b-6$.

By $k\geq3$, $s\geq k+1$, $n\geq(k-1)s+b+3$ and $n\geq2k^{2}+k+3b+3$, we deduce
\begin{align*}
&n-k(k-2)-b-3-\frac{(k-2)s+k^{2}-3k+b+4}{2(k-2)}\\
=&\frac{(2k-4)n-(k-2)s-2k^{3}+7k^{2}-(2b+11)k+3b+8}{2(k-2)}\\
\geq&\frac{(k-2)n+(k-2)((k-1)s+b+3)-(k-2)s-2k^{3}+7k^{2}-(2b+11)k+3b+8}{2(k-2)}\\
=&\frac{(k-2)n+(k-2)^{2}s-2k^{3}+7k^{2}-(b+8)k+b+2}{2(k-2)}\\
\geq&\frac{(k-2)n+(k-2)^{2}(k+1)-2k^{3}+7k^{2}-(b+8)k+b+2}{2(k-2)}\\
=&\frac{(k-2)n-k^{3}+4k^{2}-(b+8)k+b+6}{2(k-2)}\\
\geq&\frac{(k-2)(2k^{2}+k+3b+3)-k^{3}+4k^{2}-(b+8)k+b+6}{2(k-2)}\\
=&\frac{k^{3}+k^{2}+(2b-7)k-5b}{2(k-2)}\\
>&0,
\end{align*}
which leads to
\begin{align}\label{eq:4.4}
\frac{(k-2)s+k^{2}-3k+b+4}{2(k-2)}<n-k(k-2)-b-3.
\end{align}
Notice that the symmetry axis of the quadratic function $f_1(x)$ is $x=\frac{(k-2)s+k^{2}-3k+b+4}{2(k-2)}$. Combining this with \eqref{eq:4.4}, the quadratic function $f_1(x)$ is strictly increasing for
$x\geq n-k(k-2)-b-3$. Thus,
\begin{align}\label{eq:4.5}
f_1(x)\geq&f_1(n-k(k-2)-b-3)\nonumber\\
=&(k-2)(n-k(k-2)-b-3)^{2}-((k-2)s+k^{2}-3k+b+4)(n-k(k-2)-b-3)\nonumber\\
&+(k-2)(s+k)n+(b+2)n-(k^{2}-3k+2)s^{2}-(k^{3}-3k^{2}+(2b+7)k-3b-8)s\nonumber\\
&-k^{4}+3k^{3}-(2b+7)k^{2}+(3b+8)k-b^{2}-5b-6\nonumber\\
=&-(k^{2}-3k+2)s^{2}-(k^{2}+bk-b-2)s+(k-2)n^{2}\nonumber\\
&+(-2k^{3}+8k^{2}-2bk-13k+4b+10)n+k^{5}-6k^{4}+(2b+16)k^{3}\nonumber\\
&-(8b+26)k^{2}+(b^{2}+12b+24)k-2b^{2}-10b-12\nonumber\\
\geq&-(k^{2}-3k+2)\Big(\frac{n-b-3}{k-1}\Big)^{2}-(k^{2}+bk-b-2)\Big(\frac{n-b-3}{k-1}\Big)+(k-2)n^{2}\nonumber\\
&+(-2k^{3}+8k^{2}-2bk-13k+4b+10)n+k^{5}-6k^{4}+(2b+16)k^{3}\nonumber\\
&-(8b+26)k^{2}+(b^{2}+12b+24)k-2b^{2}-10b-12\nonumber\\
& \ \ \ \ \ \ \ \ \ \ (\mbox{since} \ k\geq3, \ b\geq0 \ \mbox{and} \ n\geq(k-1)s+b+3)\nonumber\\
=&\frac{1}{k-1}\Big((k-2)^{2}n^{2}-(2k^{4}-10k^{3}+(2b+22)k^{2}-(7b+29)k+7b+20)n\nonumber\\
&+k^{6}-7k^{5}+(2b+22)k^{4}-(10b+42)k^{3}+(b^{2}+21b+53)k^{2}\nonumber\\
&-(3b^{2}+25b+45)k+3b^{2}+17b+24\Big)\nonumber\\
\triangleq&\frac{1}{k-1}f_2(n).
\end{align}
Notice that
\begin{align}\label{eq:4.6}
\frac{2k^{4}-10k^{3}+(2b+22)k^{2}-(7b+29)k+7b+20}{2(k-2)^{2}}<2k^{2}+k+3b+3\leq n
\end{align}
by $k\geq3$ and $b\geq0$.

The symmetry axis of the quadratic function $f_2(n)$ is $n=\frac{2k^{4}-10k^{3}+(2b+22)k^{2}-(7b+29)k+7b+20}{2(k-2)^{2}}$. Combining this with \eqref{eq:4.6}, the quadratic function $f_2(n)$ is strictly
increasing for $n\geq2k^{2}+k+3b+3$. Hence, we obtain
\begin{align}\label{eq:4.7}
f_2(n)\geq&f_2(2k^{2}+k+3b+3)\nonumber\\
=&(k-2)^{2}(2k^{2}+k+3b+3)^{2}\nonumber\\
&-(2k^{4}-10k^{3}+(2b+22)k^{2}-(7b+29)k+7b+20)(2k^{2}+k+3b+3)\nonumber\\
&+k^{6}-7k^{5}+(2b+22)k^{4}-(10b+42)k^{3}+(b^{2}+21b+53)k^{2}\nonumber\\
&-(3b^{2}+25b+45)k+3b^{2}+17b+24\nonumber\\
=&k^{6}-k^{5}+(4b-5)k^{4}-(10b+6)k^{3}+(4b^{2}-16b+13)k^{2}\nonumber\\
&+(-18b^{2}+28b+10)k+18b^{2}+8b\nonumber\\
\triangleq&f_3(k).
\end{align}
If $k=3$, then $f_3(k)=2b+66>0$. If $k=4$, then $f_3(k)=10b^{2}+248b+1656>0$. For $k\geq5$, we possess
\begin{align*}
f_3(k)=&k^{6}-k^{5}+(4b-5)k^{4}-(10b+6)k^{3}+(4b^{2}-16b+13)k^{2}\\
&+(-18b^{2}+28b+10)k+18b^{2}+8b\\
\geq&4k^{5}+(4b-5)k^{4}-(10b+6)k^{3}+(4b^{2}-16b+13)k^{2}\\
&+(-18b^{2}+28b+10)k+18b^{2}+8b\\
\geq&(4b+15)k^{4}-(10b+6)k^{3}+(4b^{2}-16b+13)k^{2}\\
&+(-18b^{2}+28b+10)k+18b^{2}+8b\\
\geq&(10b+69)k^{3}+(4b^{2}-16b+13)k^{2}\\
&+(-18b^{2}+28b+10)k+18b^{2}+8b\\
\geq&(4b^{2}+34b+358)k^{2}+(-18b^{2}+28b+10)k+18b^{2}+8b\\
\geq&(2b^{2}+198b+1800)k+18b^{2}+8b\\
>&0.
\end{align*}
Combining these with \eqref{eq:4.3}, \eqref{eq:4.5}, \eqref{eq:4.7}, $k\geq3$ and $s\geq k+1$, we deduce
$$
\varphi_{B_2}(x)-\varphi_{B_*}(x)=(s-k)f_1(x)\geq\frac{s-k}{k-1}f_2(n)\geq\frac{s-k}{k-1}f_3(k)>0.
$$
This implies that $\varphi_{B_2}(x)>\varphi_{B_*}(x)$ for $x\geq n-k(k-2)-b-3$. Thus, we conclude
\begin{align}\label{eq:4.8}
\rho(G_2)<\rho(G_*)=\rho(K_k\vee(K_{n-k(k-1)-b-2}\cup(k(k-2)+b+2)K_1)).
\end{align}
It follows from \eqref{eq:4.1}, \eqref{eq:4.2} and \eqref{eq:4.8} that
$$
\rho(G)\leq\rho(G_1)\leq\rho(G_2)<\rho(K_k\vee(K_{n-k(k-1)-b-2}\cup(k(k-2)+b+2)K_1)),
$$
which contradicts $\rho(G)\geq\rho(K_k\vee(K_{n-k(k-1)-b-2}\cup(k(k-2)+b+2)K_1))$. This completes the proof of Theorem 1.7. \hfill $\Box$

\medskip

\section{The proof of Theorem 1.8}

\noindent{\it Proof of Theorem 1.8.} Suppose to the contrary that a connected $\frac{k}{k(k-2)+b+3}$-tough graph $G$ has no spanning tree $T$ with $\mbox{te}(T,k)\leq b$. Using Theorem 1.2, there exists
some vertex subset $S$ of $G$ such that $c(G-S)\geq(k-2)|S|+b+3$. Thus, we get $n\geq c(G-S)+|S|\geq(k-2)|S|+b+3+|S|=(k-1)|S|+b+3$. Since $G$ is connected, we conclude $S\neq\emptyset$. For convenience,
we write $|S|=s$ and $c(G-S)=c$. Then $c\geq(k-2)s+b+3$, $n\geq(k-1)s+b+3$ and $G$ is a spanning subgraph of $G_1=K_s\vee(K_{n_1}\cup K_{n_2}\cup\cdots\cup K_{n_{(k-2)s+b+3}})$ for some positive integers
$n_1,n_2,\ldots,n_{(k-2)s+b+3}$ with $\sum\limits_{i=1}^{(k-2)s+b+3}n_i+s=n$. Without loss of generality, let $n_1\geq n_2\geq\cdots\geq n_{(k-2)s+b+3}\geq1$. According to Lemma 2.4, we infer
\begin{align}\label{eq:5.1}
q(G)\leq q(G_1),
\end{align}
with equality if and only if $G=G_1$.

Let $G_2=K_s\vee(K_{n-(k-1)s-b-2}\cup((k-2)s+b+2)K_1)$. Applying Lemma 2.5, we possess
\begin{align}\label{eq:5.2}
q(G_1)\leq q(G_2),
\end{align}
with equality holding if and only if $(n_1,n_2\ldots,n_{(k-2)s+b+3})=(n-(k-1)s-b-2,1\ldots,1)$.

Recall that $|S|=s$ and $c(G-S)\geq(k-2)|S|+b+3$. It follows from $k\geq3$, $b\geq0$ and the definition of $\frac{k}{k(k-2)+b+3}$-tough graphs that
$$
\frac{k}{k(k-2)+b+3}\leq t(G)\leq\frac{|S|}{c(G-S)}\leq\frac{|S|}{(k-2)|S|+b+3},
$$
which yields $s=|S|\geq k$.

For $s=k$, we have $G_2=K_k\vee(K_{n-k(k-1)-b-2}\cup(k(k-2)+b+2)K_1)$, and so $q(G_2)=q(K_k\vee(K_{n-k(k-1)-b-2}\cup(k(k-2)+b+2)K_1))$. Together with \eqref{eq:5.1} and \eqref{eq:5.2}, we conclude
$$
q(G)\leq q(K_k\vee(K_{n-k(k-1)-b-2}\cup(k(k-2)+b+2)K_1)),
$$
with equality following if and only if $G=K_k\vee(K_{n-k(k-1)-b-2}\cup(k(k-2)+b+2)K_1)$. Then by using the condition $q(G)\geq q(K_k\vee(K_{n-k(k-1)-b-2}\cup(k(k-2)+b+2)K_1))$ of the theorem, we
obtain $G=K_k\vee(K_{n-k(k-1)-b-2}\cup(k(k-2)+b+2)K_1)$, which contradicts $G\neq K_k\vee(K_{n-k(k-1)-b-2}\cup(k(k-2)+b+2)K_1)$. Next, we shall consider $s\geq k+1$.

Recall that $G_2=K_s\vee(K_{n-(k-1)s-b-2}\cup((k-2)s+b+2)K_1)$. A direct calculation implies
\begin{align*}
2e(G_2)=&2e(K_{n-(k-2)s-b-2})+2e_{G_2}(V(K_s),V((k-2)s+b+2)K_1))\\
=&2\binom{n-(k-2)s-b-2}{2}+2s((k-2)s+b+2)\\
=&(n-(k-2)s-b-2)(n-(k-2)s-b-3)+2s((k-2)s+b+2)\\
=&k(k-2)s^{2}+(-2kn+4n+2bk+5k-2b-6)s\\
&+n^{2}-(2b+5)n+b^{2}+5b+6.
\end{align*}
Combining this with Lemma 2.6, we get
\begin{align}\label{eq:5.3}
q(G_2)\leq&\frac{2e(G_2)}{n-1}+n-2\nonumber\\
=&\frac{1}{n-1}\Big(k(k-2)s^{2}+(-2kn+4n+2bk+5k-2b-6)s\nonumber\\
&+n^{2}-(2b+5)n+b^{2}+5b+6\Big)+n-2\nonumber\\
=&\frac{1}{n-1}\Big(k(k-2)s^{2}+(-2kn+4n+2bk+5k-2b-6)s\nonumber\\
&+2n^{2}-(2b+8)n+b^{2}+5b+8\Big)\nonumber\\
\triangleq&\frac{1}{n-1}g(s).
\end{align}
Recall that $s\geq k+1$ and $n\geq(k-1)s+b+3$. Thus, we obtain $k+1\leq s\leq\frac{n-b-3}{k-1}$. It follows from $k\geq3$, $b\geq0$ and
$n\geq\frac{k^{4}+2bk^{2}+k+b^{2}+b-4}{2(k-2)}>\frac{k^{4}-2k^{3}+(b+1)k^{2}-(2b+3)k+2b+6}{(k-2)^{2}}>k^{2}+b+2$ that
\begin{align*}
g(k+1)-g\Big(\frac{n-b-3}{k-1}\Big)=&\frac{1}{(k-1)^{2}}(n-k^{2}-b-2)\Big((k-2)^{2}n-k^{4}\\
&+2k^{3}-(b+1)k^{2}+(2b+3)k-2b-6\Big)\\
>&0,
\end{align*}
which implies that $g(s)$ attains its maximum value at $s=k+1$ for $k+1\leq s\leq\frac{n-b-3}{k-1}$. Together with \eqref{eq:5.3}, $k\geq3$, $b\geq0$ and $n\geq\frac{k^{4}+2bk^{2}+k+b^{2}+b-4}{2(k-2)}$,
we infer
\begin{align}\label{eq:5.4}
q(G_2)\leq&\frac{1}{n-1}g(s)\nonumber\\
\leq&\frac{1}{n-1}g(k+1)\nonumber\\
=&\frac{1}{n-1}\Big(k(k-2)(k+1)^{2}+(-2kn+4n+2bk+5k-2b-6)(k+1)\nonumber\\
&+2n^{2}-(2b+8)n+b^{2}+5b+8\Big)\nonumber\\
=&\frac{1}{n-1}\Big(2n^{2}-(2k^{2}-2k+2b+4)n+k^{4}+(2b+2)k^{2}-3k+b^{2}+3b+2\Big)\nonumber\\
=&2(n-k(k-2)-b-3)-\frac{2(k-2)n-k^{4}-2bk^{2}-k-b^{2}-b+4}{n-1}\nonumber\\
\leq&2(n-k(k-2)-b-3).
\end{align}

Note that $K_k\vee(K_{n-k(k-1)-b-2}\cup(k(k-2)+b+2)K_1)$ contains $K_{n-k(k-2)-b-2}$ as its proper subgraph. In terms of Lemma 2.4, we conclude
\begin{align}\label{eq:5.5}
q(K_k\vee(K_{n-k(k-1)-b-2}\cup(k(k-2)+b+2)K_1))>&q(K_{n-k(k-2)-b-2})\nonumber\\
=&2(n-k(k-2)-b-3).
\end{align}
It follows from \eqref{eq:5.1}, \eqref{eq:5.2}, \eqref{eq:5.4} and \eqref{eq:5.5} that
$$
q(G)\leq q(G_1)\leq q(G_2)<q(K_k\vee(K_{n-k(k-1)-b-2}\cup(k(k-2)+b+2)K_1)),
$$
which contradicts $q(G)\geq q(K_k\vee(K_{n-k(k-1)-b-2}\cup(k(k-2)+b+2)K_1))$. Theorem 1.8 is proved. \hfill $\Box$

\section*{Declaration of competing interest}

The author declares that he has no known competing financial interests or personal relationships that could have appeared to influence the work reported in this paper.

\section*{Data availability}

No data was used for the research described in the article.

\section*{AI Contribution}

No use of AI tools were made in the preparation of this article.



\section*{Acknowledgments}

This work was supported by the Natural Science Foundation of Jiangsu Province (Grant No. BK20241949).


\begin{thebibliography}{9999}

\bibitem {B} V. Benediktovich, Spectral conditions of pancyclicity for $t$-tough graphs, Discrete Appl. Math. 365 (2025) 130--137.

\bibitem {Os} S. O, Spectral radius and matchings in graphs, Linear Algebra Appl. 614 (2021) 316--324.

\bibitem {Ws2} J. Wu, Sufficient conditions for a graph with minimum degree to have a component factor, Proc. Rom. Acad. Ser. A Math. Phys. Tech. Sci. Inf. Sci. 27(1) (2026) 3--10.

\bibitem {Zs1} S. Zhou, Sufficient conditions for a graph with minimum degree to be $k$-critical with respect to an odd $[1,b]$-factor, Discrete Math. 349(11) (2026) 115251.

\bibitem {Zs2} S. Zhou, Some spectral conditions for star-factors in bipartite graphs, Discrete Appl. Math. 369 (2025) 124--130.

\bibitem {Zs3} S. Zhou, Spanning subgraphs and spectral radius in graphs, Aequationes Math. 100(1) (2026) 1.

\bibitem {PK} S. Pirzada, S. Khan, On Zagreb index, signless Laplacian eigenvalues and signless Laplacian energy of a graph, Comput. Appl. Math. 42(4) (2023) 152.

\bibitem {Ws1} J. Wu, Some results on the $k$-strong parity property in a graph, Comput. Appl. Math. 45(4) (2026) 138.

\bibitem {ZL} S. Zhou, H. Liu, Perfect matchings and $A_{\alpha}$-spectral radius in 1-binding graphs, Discrete Appl. Math. 397 (2027) 35--42.

\bibitem {ZLB} S. Zhou, H. Liu, Q. Bian, $A_{\alpha}$-spectral radius and path-factor covered graphs, AIMS Math. 11(9) (2026) 28646--28664.

\bibitem {ZZZL} S. Zhou, Y. Zhang, T. Zhang, H. Liu, Toughness and $A_{\alpha}$-spectral radius in graphs, Filomat 40(5) (2026) 1883--1892.

\bibitem {KKS} A. Kaneko, M. Kano, K. Suzuki, Spanning trees with leaf distance at least four, J. Graph Theory 55 (2007) 83--90.

\bibitem {Wc} J. Wu, Characterizing spanning trees via the size or the spectral radius of graphs, Aequationes Math. 98(6) (2024) 1441--1455.

\bibitem {WZL} J. Wu, S. Zhou, H. Liu, A spectral condition for spanning trees with restricted degrees in bipartite graphs, Proc. Rom. Acad. Ser. A Math. Phys. Tech. Sci. Inf. Sci. 27(1) (2026) 19--24.

\bibitem {OS} K. Ota, T. Sugiyama, Forbidden subgraphs and the existence of spanning $k$-trees, Discrete Math. 310 (2010) 3506--3511.

\bibitem {FMMMTY} M. Furuya, S. Maezawa, R. Matsubara, H. Matsuda, S. Tsuchiya, T. Yashima, Degree sum conditions for the existence of spanning $k$-trees in star-free graphs, Discuss. Math. Graph Theory 42(1) (2022) 5--13.

\bibitem {R} D. Rose, On simple characterizations of $k$-trees, Discrete Math. 7 (1974) 317--322.

\bibitem {GL} X. Gu, M. Liu, A tight lower bound on the matching number of graphs via Laplacian eigenvalues, European J. Combin. 101 (2022) 103468.

\bibitem {FGHL} D. Fan, S. Goryainov, X. Huang, H. Lin, The spanning $k$-trees, perfect matchings and spectral radius of graphs, Linear Multilinear Algebra 70 (2022) 7264--7275.

\bibitem {ZW} S. Zhou, J. Wu, Spanning $k$-trees and distance spectral radius in graphs, J. Supercomput. 80(16) (2024) 23357--23366.

\bibitem {WZ} J. Wu, S. Zhou, Spanning $k$-trees, odd $[1,b]$-factors and spectral radius in binding graphs, https://arxiv.org/pdf/2507.07301

\bibitem {LFS} R. Liu, A. Fan, J. Shu, Spectral extremal problems on factors in tough graphs, and beyond, Discrete Math. 348 (2025) 114593.

\bibitem {OOO} Y. Ohnishi, K. Ota, K. Ozeki, A note on total excess of spanning trees, AKCE Int. J. Graphs Comb. 8(1) (2011) 97--103.

\bibitem {MTY} S. Maezawa, M. Tsugaki, T. Yashima, Closure and spanning trees with bounded total excess, Graphs Combin. 37 (2021) 805--822.

\bibitem {Oa} K. Ozeki, A toughness condition for a spanning tree with bounded total excesses, Graphs Combin. 31 (2015) 1679--1688.

\bibitem {Win} S. Win, On a connection between the existence of $k$-trees and the toughness of a graph, Graphs Combin. 5 (1989) 201--205.

\bibitem {EZ} M. Ellingham, X. Zha, Toughness, trees, and walks. J Graph Theory 33(3) (2000) 125--137.

\bibitem {EOO} H. Enomoto, Y. Ohnishi, K. Ota, Spanning trees with bounded total excess, Ars Combin. 102 (2011) 289--295.

\bibitem {Za} S. Zhou, A result on spanning trees with bounded total excess, Discrete Appl. Math. 388 (2026) 130--135.

\bibitem {GZ} Y. Gong, P. Zhang, Spectral and size conditions for the existence of spanning trees with bounded total excess, Proc. Rom. Acad. Ser. A Math. Phys. Tech. Sci. Inf. Sci. 27(3) (2026) 245--252.

\bibitem {ZLLW} L. Zheng, S. Li, X. Luo, G. Wang, Some sufficient conditions for a graph with minimum degree to be $k$-factor-critical, Discrete Appl. Math. 348 (2024) 279--291.

\bibitem {LF} Q. Li, K. Feng, On the largest eigenvalue of a graph, Acta Math Appl Sin Chinese Ser 2 (1979) 167--175.

\bibitem {SYZL} Y. Shen, L. You, M. Zhang, S. Li, On a conjecture for the signless Laplacian spectral radius of cacti with given matching number, Linear Multilinear Algebra 65 (2017) 457--474.

\bibitem {Das} K. Das, Maximizing the sum of the squares of the degrees of a graph, Discrete Math. 285 (2004) 57--66.

\bibitem {YYSX} L. You, M. Yang, W. So, W. Xi, On the spectrum of an equitable quotient matrix and its application, Linear Algebra Appl. 577 (2019) 21--40.



\end{thebibliography}
\end{document}